\documentclass[1p]{elsarticle}

\usepackage{lineno,hyperref}
\modulolinenumbers[1]

\journal{Applications of Fractional Calculus to Modeling in Dynamics and Chaos}

\usepackage{epstopdf}
\usepackage{amsmath,amsfonts,amssymb}

\usepackage[section]{placeins}
\usepackage{multirow}
\usepackage{mathtools}

\newenvironment{Proof}[1]{
\:{\textbf{Proof.}}\:\:{\rmfamily #1}
}{
\hfill$\blacksquare$
}

\newtheorem{definition}{Definition}[section]
\newtheorem{thm}{Theorem}[section]
\newtheorem{rem}{Remark}[section]
\newtheorem{lemma}{Lemma}[section]
\newtheorem{cor}{Corollary}[section]

\usepackage{mathtools}

\providecommand{\norma}[1]{\Bigg\lVert#1 \Bigg\rVert}

\providecommand{\abs}[1]{\lvert#1\rvert}
\providecommand{\norm}[1]{\lVert#1\rVert}

\providecommand{\norma}[1]{\Bigg\lVert#1 \Bigg\rVert}

\providecommand{\abs}[1]{\lvert#1\rvert}
\providecommand{\norm}[1]{\lVert#1\rVert}

\usepackage{hyperref}

\begin{document}

\begin{frontmatter}
\title{Stability of Fractional Nonlinear Systems with Mittag-Leffler Kernel and Design of State Observers}

\author[Ana]{O. Mart\'inez-Fuentes \corref{cor1}}
\ead{oscar.martinezfu@anahuac.mx}
\author[Iber]{Sergio M. Delf\'in-Prieto}
\ead{sergio.delfin@itnl.edu.mx}

\cortext[cor1]{Corresponding author}
\address[Ana]{School of Engineering, Universidad An\'ahuac-Xalapa. Circuito Arco Sur s/n, Col. Lomas Verdes, C.P. 91098, Xalapa, Veracruz, Mexico}
\address[Iber]{Instituto Tecnol\'ogico de Nuevo Le\'on. Av. Eloy Cavazos 2001, Tolteca 67170. Guadalupe, N. L, Mexico.}

\begin{abstract}
Atangana and Baleanu proposed a new fractional derivative with non-local and no-singular Mittag–Leffler kernel to solve some problems proposed by researchers in the field of fractional calculus. This new derivative is better to describe essential aspects of non-local dynamical systems. We present some results regarding Lyapunov stability theory, particularly the Lyapunov Direct Method for fractional-order systems modeled with Atangana-Baleanu derivatives and some significant inequalities that help to develop the theoretical analysis. As applications in control theory, some algorithms of state estimation are proposed for linear and nonlinear fractional-order systems.
\end{abstract}

\end{frontmatter}

\section{Introduction}\label{intro}

The fractional calculus is the study generalized differentiation and integration of arbitrary real or complex order or even variable. This branch of mathematics is referred to as differintegration. It was only at the beginning of the XIX century that the generalization of differentiation to real or complex orders was formalized with the work of Liouville, Riemann and Letnikov \cite{oldham1974fractional}. Oldham and Spanier wrote the first work devoted exclusively to the subject of fractional calculus \cite{oldham1974fractional}.

There are many definitions of fractional-order integrals and derivatives. By starting with the most often used: the Riemann–Liouville integral. The fractional Riemann–Liouville integrals are defined as the general case of iterated integrals of integer order to an arbitrary positive number. Meanwhile, the fractional Riemann–Liouville derivatives are determined to be derivatives (in the classical sense) of fractional integrals. Another frequently used definition of fractional derivative is the Caputo definition. In this case, fractional integrals are again defined using the definition of Riemann-Liouville, while fractional derivatives are determined to be fractional integrals of classic derivatives.

The applications of fractional calculus have applied in the engineering field and physics, for example, the viscoelasticity phenomena. The complexity of viscoelastic materials that occurs in the linear domain was explained by the influence of memory on the dynamics and the possibility of modeling these effects using the fractional calculus \cite{koeller1984applications,bagley1983theoretical}.

The heat conduction problem and one proposed solution by means of the heat equation with memory are established under some general and reasonable conditions in \cite{yong2005heat}. The double Laplace transform method solves a fractional heat equation subject to certain initial and boundary conditions \cite{anwar2013fractional}. The second law of thermodynamics ensures that the heat flows from hot to cold regions, and this condition is analyzed in the context of the fractional calculus in \cite{vazquez2011fractional}.

Some potential applications and theoretical implications for the feedback control theory with fractional differential equations have been done in different approaches. The internal and external stability properties of finite-dimensional linear fractional-order systems are studied in \cite{matignon1996stability}, and the controllability and observability are revisited in \cite{matignon1996some}.

The study of controllability and the null controllability of fractional dynamical systems with limited control is established in \cite{balachandran2016note,nirmala2017null}. On the other hand, the formulation of optimal control is carried out in \cite{tricaud2010approximate,djennoune2013optimal,trujillo2018optimal,wei2017fractional}.

In \cite{gomez2014physical} and \cite{guia2013analysis}, the authors emphasize that fractional differentiation with respect to time can be interpreted as the existence of memory effects that correspond to intrinsic dissipation in a system. This approach can be observed in the time constant and transitory response of a capacitor charge and discharge. In \cite{jesus2012application}, the use of integer and fractional electrical elements is described for modeling two electrochemical systems. The first system consists of botanical elements, while the second system is implemented with fractal electrodes through electrolyte processes.

Other areas where the theory of fractional differential equations has been applied are the mechanics and the fractal analysis. The generalization of the Newtonian equation is derived from the memory of the fractional case \cite{baleanu2010newtonian}. The Euler-Lagrange equations and the methods of generalized mechanics are applied to obtain the Hamilton equations for free motion from fractional-order calculus perspective \cite{klimek2002lagrangean,baleanu2005lagrangian}.

The development of fractional-derivative mechanics by deriving a modified Hamilton principle and the Hamilton-Jacobi equation using generalized mechanics with fractional and higher-order derivatives \cite{riewe1997mechanics,riewe1996nonconservative,gomez2015modeling}. The fractal calculus involving the non-local derivatives is generalized, and the scaling properties of the local and non-local derivatives are studied in \cite{golmankhaneh2016non}.

The physical interpretation of fractional calculus is an open issue. However, the convolution kernels in physical laws are of great importance in many models, and they can be useful in the physical interpretation. This idea can be justified due to the kernels govern different physical phenomena, such as memory effects. For example, the role of memory functions for compliance retardation and modulus relaxation in viscoelastic materials is examined. Also, the mathematical tools of fractional calculus lead to results that demonstrate the fractional behavior of memory properties of some signals \cite{machado2011and}.

Based on the scale invariance property, which exists for fractals with discrete structure, it becomes possible to understand the geometrical and physical meaning of the fractional derivatives and integrals with complex fractional exponents. The correct form of this complex structure, which can enter into kinetic equation with fractional derivative, has been found in \cite{nigmatullin2005there}.

The solutions of a differential equation in the sense of Caputo incorporate and describe long-term memory effects (attenuation or dissipation). These effects are linked to the Mittag–Leffler function. If the fractional-order is less than one, the fractional differentiation with respect to time represents the non-local displacement effect of energy dissipation (internal friction) \cite{gomez2015modeling}. In \cite{heymans2006physical} demonstrated on a series of examples that it is possible to attribute physical meaning to initial conditions expressed in terms of Riemann-Liouville fractional derivatives.

In the case of the fractional-order models of viscoelasticity (the Zener model), the initial conditions are zero. In such cases, the use of the Riemann-Liouville derivatives, the Gr\"{u}nwald-Letnikov derivative, and the Caputo derivative are equivalent.

These approaches of differintegration prevail with fruitful applications in engineering until some research realized that some fractional derivatives present a singular kernel in their basic definition. The Caputo definition and other fractional derivatives own singular kernel, which means that mathematical models involving these sorts of fractional derivatives do not consider the complete memory \cite{mekkaoui2019new}. That is, the ordinary derivative is a local operator, and dynamic processes are memoryless. In contrast, the fractional operators defined are non-local. Therefore the dynamic processes modeled with those operators hold a degree of memory \cite{tarasov2018no,yepez2018numerical,abdeljawad2018fractional}.

Atangana and Baleanu suggested a fractional derivative involving Mittag-Leffler function in its kernel, and such derivative is now called Atangana-Baleanu fractional derivative \cite{atangana2016new,baleanu2018some}. This operator aimed to overcome the problem of the singularity of the kernel and involve full memory effect into the fractional derivative \cite{baleanu2018nonlinear}. This feature was introduced first by Caputo and Fabrizio by adding an exponential function in its kernel \cite{caputo2016applications,li2019stability,hristov2019linear}.

In this new perspective, the classical fractional calculus and the Mittag-Leffler non-singular kernel approach, both complement each other in making an effort to describe and highlight the covert properties of non-local dynamical processes \cite{srivastava2009fractional,gomez2018fractional,owolabi2019modelling,kumar2019analysis,avalos2019dynamics}.

This chapter discusses the interaction between two research fields: the fractional-order dynamic systems described by derivative with Mittag-Leffler kernel (FOSMLK) and the design of nonlinear observer for a class of nonlinear systems. On the one hand, FOSMLK refers to a class of nonlinear systems that are defined with a fractional derivative with Mittag-Leffler kernel, which brings out more clearly the non-locality of fractional derivatives. On the other hand, we revisit the concept of stability in the sense of Lyapunov with the Atangana-Baleanu fractional derivative \cite{abdeljawad2017lyapunov}. With the obtaining results, we addressed the design of observers for a class of nonlinear systems.

The question of estimating states for classical integer-order dynamic systems has been extensively studied in the literature. Observers are regarded as essential algorithms in control theory, with the principal target of reconstructing the states of a dynamic process. Note that these algorithms have been widely used in other advanced tasks, such as control problems. The pioneer contribution of Luenberger breaks out the career to overcome different scenarios in the designing observers for linear and nonlinear systems \cite{luenberger1971introduction,ciccarella1993luenberger}. This approach is based on the contribution of Kalman \cite{kalman1960contributions}, which is related to the structural properties of observability and controllability for linear systems. However, in the case of observers for nonlinear systems, a well-known procedure is widely used: the Lyapunov analysis.

One can make a first observation about the work of Lyapunov, i.e., The General Problem of Motion Stability \cite{lyapunov1992general} that encompasses two methods for stability analysis published in 1892. The first method, also named as linearization method, brings out conclusions about a nonlinear system's local stability around an equilibrium point from the stability properties of its linear approximation. The direct method, nowadays known as the Lyapunov direct method, is not restricted to local motion and determines the stability properties of a nonlinear system by constructing a scalar energy-like function for the system and examining the function's time variation.

The procedure of Lyapunov analysis is constructing a quadratic-like function and then analyze the derivative. Roughly speaking, strict Lyapunov functions are characterized by having negative definite time derivatives along all trajectories of the system; meanwhile, non-strict Lyapunov functions have negative semidefinite derivatives along the trajectories. Together, the first method and the direct method constitute the so-called Lyapunov stability theory.

Now, the generalization of actual results to the fractional-order framework represents an opportunity area of research.  To the best of knowledge, the general inquiry of designing observers for fractional-order systems with Mittag-Leffler kernel has not been investigated either a constructive method for the analysis stability of fractional differential equations with the Atangana-Baleanu derivative.

The objective of this chapter is to present the Lyapunov stability theory for a class of nonlinear systems and illustrate its use in the analysis and the design of observers for nonlinear systems modeled by the Atangana-Baleanu fractional derivative. Some results about the Mittag-Leffler stability are given in terms of nonautonomous time-invariant systems. Lyapunov-like inequalities are obtained and employed to conclude the convergence of a dynamic system. Also, we introduce a family of observers, which we call them Mittag-Leffler observers. These sorts of observers are characterized by the convergence rate and are bounded by a Mittag-Leffler function, which is a suitable behavior for fractional-order dynamic systems.

\section{Fractional Derivatives with Mittag-Leffler Kernel}

In this section, we present the definitions of the fractional derivatives with Mittag-Leffler kernel introduced by Atangana and Baleanu in \cite{atangana2016new}, the classical Riemann-Liouville integral and some useful results in the solution of fractional differential equations.

\begin{definition}\label{int_RL}{\rm
Let $\alpha>0$. The operator $\prescript{RL}{t_{0}}{I}^{\alpha}_{t}$ defined on $L^{1}[t_0,T]$ by
\begin{equation}\label{Integral_RIEMLIOU}
\prescript{RL}{t_{0}}{I}^{\alpha}_{t} f(t)=\frac{1}{\Gamma(\alpha)}\int_{t_0}^{t}\left(t-\tau \right)^{\alpha-1}f(\tau)\:d\tau
\end{equation}
for $t_0\leq t\leq T$, is called the Riemann-Liouville fractional integral operator of order $\alpha$.
\hfill{$\Box$}
}\end{definition}

\begin{rem}
In Definition \ref{int_RL}, $\Gamma(z)$ represents the classical gamma function
\begin{equation}\label{gamma}
\Gamma(z)=\int_{0}^{\infty}t^{z-1}e^{-t}\:dt.
\end{equation}
It is essential to mention that $\operatorname{Re}(z)\in\mathbb{R}\setminus\left(\{0\}\cup\mathbb{Z}^{-}\right)$.
\end{rem}

\begin{definition}\label{1def:ABR}{\rm
The Atangana-Baleanu (AB) fractional derivative of $f(t)$ in the Riemann-Liouville (R) sense is defined as follows
\begin{equation}\label{ABR_deriv}
\prescript{ABR}{t_{0}}{\mathcal{D}}^{\alpha}_{t} f(t)=\frac{B(\alpha)}{1-\alpha}\frac{d}{dt}\int_{t_0}^{t} f(\tau) E_{\alpha}\left( \dfrac{-\alpha}{1-\alpha}\left(t-\tau\right)^{\alpha} \right)\:d\tau
\end{equation}
for $0<\alpha<1$, $t_0<t<T$, and $f\in {L}^{1}(t_0,T)$.
\hfill{$\Box$}
}\end{definition}

\begin{definition}\label{1def:ABC}{\rm
The AB fractional derivative of $f(t)$ in the Caputo (C) sense is defined by
\begin{equation}\label{ABC_deriv}
\prescript{ABC}{t_{0}}{\mathcal{D}}^{\alpha}_{t} f(t)=\frac{B(\alpha)}{1-\alpha}\int_{t_0}^{t} f^\prime(\tau) E_{\alpha}\left( \dfrac{-\alpha}{1-\alpha}\left(t-\tau\right)^{\alpha} \right)\:d\tau
\end{equation}
for $0<\alpha<1$, $t_0<t<T$, and $f$ a differentiable function on $[t_0,T]$ such that $f^\prime\in {L}^{1}(t_0,T)$.
\hfill{$\Box$}
}\end{definition}

In the Definitions \ref{1def:ABR} and \ref{1def:ABC}, $B(\alpha)$ denotes a real-valued normalization function satisfying $B(\alpha)>0$ and $B(0)=B(1)=1$. Furthermore, $E_\alpha(\cdot)$ is the one parameter Mittag-Leffler function defined by
\begin{equation}\label{ML_function}
E_{\alpha}(z)=\sum_{k=0}^{\infty}\frac{z^k}{\Gamma\left(\alpha k+1 \right)},
\end{equation}
whenever the series converges \cite{gorenflo2014mittag}. The function \eqref{ML_function} is a particular case of the two-parameter Mittag–Leffler function given by
\begin{equation}\label{ML_function}
E_{\alpha,\beta }(z)=\sum_{k=0}^{\infty}\frac{z^k}{\Gamma\left(\alpha k+\beta \right)}.
\end{equation}
In both cases, the parameters $\alpha$ and $\beta$ are such that $\operatorname{Re}(\alpha)>0$, $\operatorname{Re}(\beta)>0$. In addition to the previous definitions of fractional-derivatives, the AB fractional integral operator is well defined, and some results can be obtained together with the AB derivative.

\begin{definition}\label{int_ABC}{\rm
The AB fractional integral of $f(t)$ is defined by
\begin{equation}\label{AB_integral}
\prescript{AB}{t_{0}}{I}^{\alpha}_{t} f(t)=\frac{1-\alpha}{B\left(\alpha\right)}f(t)+\frac{\alpha}{B\left(\alpha\right)} \prescript{RL}{t_{0}}{I}^{\alpha}_{t} f(t)
\end{equation}
for $\alpha>0$ and $\prescript{RL}{t_{0}}{I}^{\alpha}_{t}$ denotes the Riemann-Liouville integral.
\hfill{$\Box$}
}\end{definition}

\begin{cor}\cite{baleanu2018some}\label{corolario_integral_deri}
The AB fractional integral and the ABC fractional derivative satisfy the following Newton-Leibniz formula:
\begin{equation}\label{New_Lei}
\prescript{AB}{t_{0}}{I}^{\alpha}_{t}\left(\prescript{ABC}{t_{0}}{\mathcal{D}}^{\alpha}_{t}f(t) \right)=f(t)-f(t_0)
\end{equation}
\end{cor}

A fundamental tool for solving differential equations is the Laplace Transform $F(s)=\mathcal{L}\left\{f(t)\right\}=\displaystyle \int_{0}^{\infty}e^{-st}f(t)\:dt$ for functions $f$ with exponential order. The Laplace transform for ABC and  ABR derivatives, as well as for the integral AB can be derived by simple calculations.

\begin{equation}\label{Lapla_ABC}
\mathcal{L} \left\{ \prescript{ABC}{{0}}{\mathcal{D}}^{\alpha}_{t} f(t)\right\}=\frac{B(\alpha)}{1-\alpha}\dfrac{s^{\alpha-1}}{s^{\alpha}+\dfrac{\alpha}{1-\alpha}}\left[sF(s)-f(0)\right]
\end{equation}
\begin{equation}\label{Lapla_ABR}
\mathcal{L} \left\{ \prescript{ABR}{{0}}{\mathcal{D}}^{\alpha}_{t} f(t)\right\}=\frac{B(\alpha)}{1-\alpha}\dfrac{s^{\alpha}}{s^{\alpha}+\dfrac{\alpha}{1-\alpha}}F(s)
\end{equation}
\begin{equation}\label{Lapla_integral}
\mathcal{L} \left\{\prescript{AB}{t_{0}}{I}^{\alpha}_{t} f(t)\right\}=\frac{1}{B(\alpha)}\left[\frac{(1-\alpha)s^{\alpha}+\alpha}{s^{\alpha}} \right]F(s)
\end{equation}

As we have seen so far, the Mittag-Leffler function is quite relevant in fractional differential equations. When solving differential equations using the Laplace transform method, the following function is handy. Let $\alpha>0,\beta>0$ and $\lambda\in \mathbb{C}$, then
\begin{equation}\label{Lapla_ML}
\mathcal{L} \left\{t^{\beta-1}E_{\alpha,\beta}\left(-\lambda t^\alpha\right)\right\}=\frac{s^{\alpha-\beta}}{s^\alpha+\lambda}
\end{equation}

As we know, the classical Riemann-Liouville derivative can be written as a function of the classical Caputo derivative. ABC and ABR derivatives can also be related, as we will see below. This equation allows obtaining results for one derivative and extending them to the other with particular care and conditions, as we will see later.

\begin{lemma}\label{lemma_1_important}
\begin{equation}\label{relation_derivatives}
\prescript{ABR}{t_{0}}{\mathcal{D}}^{\alpha}_{t} f(t)=\prescript{ABC}{t_{0}}{\mathcal{D}}^{\alpha}_{t} f(t)+\frac{B(\alpha)}{1-\alpha}f(0)E_{\alpha}\left( \dfrac{-\alpha}{1-\alpha} t^\alpha \right)
\end{equation}

\begin{Proof}
{\rm From the equation \eqref{ABR_deriv} and the convolution properties, one has that the ABR derivative can be written as follows:
\begin{equation}\label{}
 \prescript{ABR}{t_{0}}{\mathcal{D}}^{\alpha}_{t} f(t)=\frac{B(\alpha)}{1-\alpha}\frac{d}{dt}\int_{t_0}^{t} f(t-\tau) E_{\alpha}\left( \dfrac{-\alpha}{1-\alpha}\tau^{\alpha} \right)\:d\tau.
\end{equation}
By using the Leibniz's rule for differentiation under the integral sign \cite{flanders}, it follows that
\begin{eqnarray*}
 \prescript{ABR}{t_{0}}{\mathcal{D}}^{\alpha}_{t} f(t)&=&\frac{B(\alpha)}{1-\alpha}\left[\int_{t_0}^{t} f^\prime(t-\tau) E_{\alpha}\left( \dfrac{-\alpha}{1-\alpha}\tau^{\alpha} \right)\:d\tau+f(0)E_{\alpha}\left( \dfrac{-\alpha}{1-\alpha} t^\alpha \right)\right]\\
 &=&\frac{B(\alpha)}{1-\alpha}\int_{t_0}^{t} f^\prime(\tau) E_{\alpha}\left( \dfrac{-\alpha}{1-\alpha}\left(t-\tau\right)^{\alpha} \right)\:d\tau+  \frac{B(\alpha)}{1-\alpha}f(0)E_{\alpha}\left( \dfrac{-\alpha}{1-\alpha} t^\alpha \right).
\end{eqnarray*}
This concludes the proof considering that the first term is the ABC derivative and an extra term as we expected.
}
\end{Proof}
\end{lemma}

We conclude this section with one inequality based on the expression \eqref{relation_derivatives} to determine the comparison between fractional derivatives.

\begin{lemma}\label{lema_ineq}
Let $f$ a function that satisfies the conditions for existence for ABR and ABC derivatives, and such that $f(0)\geq 0$, then
\begin{equation}\label{ineq_deriv}
\prescript{ABC}{t_{0}}{\mathcal{D}}^{\alpha}_{t} f(t)\leq \prescript{ABR}{t_{0}}{\mathcal{D}}^{\alpha}_{t} f(t)
\end{equation}

\begin{Proof}
{\rm
From the Lemma \ref{lemma_1_important}, we have
$$
\prescript{ABC}{t_{0}}{\mathcal{D}}^{\alpha}_{t} f(t)=\prescript{ABR}{t_{0}}{\mathcal{D}}^{\alpha}_{t} f(t)- \frac{B(\alpha)}{1-\alpha}f(0)E_{\alpha}\left( \dfrac{-\alpha}{1-\alpha} t^\alpha \right)
$$
Considering the positivity of the Mittag-Leffler function for $0<\alpha<1$ and $f(0)\geq 0$, we get the desired result.
}\end{Proof}
\end{lemma}

\section{On Fractional-Order Nonlinear Systems and Mittag-Leffler Stability}

In this section, we discuss some results about the fractional nonautonomous system
\begin{equation}\label{system_non}
\prescript{ABC}{t_{0}}{\mathcal{D}}^{\alpha}_{t} x(t)=f(t,x_t)
\end{equation}
with initial condition $x(t_0)$, $0<\alpha<1$, $f:[t_0,\infty]\times\Omega\to \mathbb{R}^n$ is piecewise continuous in $t$ and locally Lipschitz in $x$ on $[t_0,\infty]\times \Omega$, and $\Omega\in \mathbb{R}^n$ is a domain that contains the origin $x=0$.

\begin{rem}
Some results about uniqueness and existence of the solution to the fractional-order system \eqref{system_non} are summarized in \cite{atangana2016new,baleanu2018nonlinear,owolabi2019modelling}.
\end{rem}

\begin{lemma}\label{lemma_bound_integralAB}
For the real-valued continuous $f(t,x_t)$ in \eqref{system_non}, we have
\begin{equation}\label{norma_IAB}
\norm{\prescript{AB}{t_{0}}{I}^{\alpha}_{t} f(t,x_t)}\leq \prescript{AB}{t_{0}}{I}^{\alpha}_{t} \norm{f(t,x_t)}
\end{equation}
where $\alpha>0$.

\begin{Proof}
{\rm
From the integral \eqref{AB_integral} and the properties of the norm, we have that
\begin{eqnarray*}
\norm{\prescript{AB}{t_{0}}{I}^{\alpha}_{t} f(t,x_t)}&=&\norma{\frac{1-\alpha}{B\left(\alpha\right)}f(t)+\frac{\alpha}{B\left(\alpha\right)} \prescript{RL}{t_{0}}{I}^{\alpha}_{t} f(t)}\\
&\leq & \frac{1-\alpha}{B\left(\alpha\right)}\norm{f(t)}+\frac{\alpha}{B\left(\alpha\right)} \norm{\prescript{RL}{t_{0}}{I}^{\alpha}_{t} f(t)}\\
&\leq & \frac{1-\alpha}{B\left(\alpha\right)}\norm{f(t)}+\frac{\alpha}{B\left(\alpha\right)} \prescript{RL}{t_{0}}{I}^{\alpha}_{t} \norm{f(t)}\\
&=&\prescript{AB}{t_{0}}{I}^{\alpha}_{t} \norm{f(t,x_t)}
\end{eqnarray*}

The last inequality is valid because the Riemann-Liouville integral is bounded (see Lemma 2 of \cite{li2010stability}).

}
\end{Proof}
\end{lemma}

\begin{thm}
If $x=0$ is an equilibrium point of system \eqref{system_non}, $f$ is Lipschitz on $x$ with Lipschitz constant $\kappa$ and is piecewise continuous with respect to $t$, then the solution of \eqref{system_non} satisfies
\begin{equation}\label{bound_1}
\norm{x(t)}\leq \frac{B(\alpha)}{B(\alpha)-\kappa(1-\alpha)}\norm{x(t_0)}E_{\alpha} \left(\frac{\kappa \alpha}{B(\alpha)-\kappa(1-\alpha)}      t^\alpha\right),
\end{equation}
where $0<\alpha<1$.

\begin{Proof}{\rm
By applying the AB integral \eqref{AB_integral} to both sides of \eqref{system_non}, it follows from Corollary \ref{corolario_integral_deri}, Lemma \ref{lemma_bound_integralAB} and the Lipschitz condition that
\begin{equation*}
\begin{aligned}
\abs{\norm{x(t)}-\norm{x(t_0)}}\leq \norm{x(t)-x(t_0)}=\norm{\prescript{AB}{t_{0}}{I}^{\alpha}_{t} f(t,x_t)}\\
\leq \prescript{AB}{t_{0}}{I}^{\alpha}_{t} \norm{f(t,x_t)}\leq \kappa  \prescript{AB}{t_{0}}{I}^{\alpha}_{t} \norm{x(t)}.
\end{aligned}
\end{equation*}

Let $M(t)$ a non-negative function, then
\begin{equation}\label{eq_auxM}
\norm{x(t)}-\norm{x(t_0)}+M(t)=\kappa  \prescript{AB}{t_{0}}{I}^{\alpha}_{t} \norm{x(t)}.
\end{equation}
By applying the Laplace transform to \eqref{eq_auxM}, it follows that
\begin{equation}\label{new_eq}
\norm{X(s)}\left[\lambda_1s^{\alpha}-\lambda_2\right]=\norm{x(t_0)}s^{\alpha-1}-s^\alpha M(s)
\end{equation}
where
$$
\lambda_1=1-\frac{\kappa}{B(\alpha)}(1-\alpha),\qquad \lambda_2=\frac{\kappa \alpha}{B(\alpha)},
$$
$$
\norm{X(s)}=\int_{0}^{\infty}e^{-st}\norm{x(t)}\:dt,\: M(s)=\int_{0}^{\infty}e^{-st}m(t)\:dt.
$$

Hence, from the equation \eqref{new_eq}, one has that:
\begin{equation}\label{laplace_thm}
\norm{X(s)}=\frac{\norm{x(t_0)}}{\lambda_1}\frac{s^{\alpha-1}}{s^\alpha-\lambda_1^{-1}\lambda_2}-\frac{\lambda_2}{\lambda^2_1}\frac{1}{s^\alpha-\lambda_1^{-1}\lambda_2}M(s)-\frac{1}{\lambda_1}M(s)
\end{equation}

Applying the inverse Laplace transform to \eqref{laplace_thm} yields
\begin{equation}
\norm{x(t)}=\frac{\norm{x(t_0)}}{\lambda_1} E_{\alpha} \left(\frac{\lambda_2}{\lambda_1}t^{\alpha}   \right)-\frac{\lambda_2}{\lambda^2_1}t^{\alpha-1}E_{\alpha,\alpha} \left(\frac{\lambda_2}{\lambda_1}t^{\alpha}   \right)\ast M(t)-\frac{1}{\lambda_1}M(t)
\end{equation}
where $\ast$ is convolution operator. Finally, since both $t^{\alpha-1}$ and $E_{\alpha,\alpha} \left(\frac{\lambda_2}{\lambda_1}t^{\alpha}   \right)$ are nonnegative functions, we conclude that
$$
\norm{x(t)}\leq \frac{B(\alpha)}{B(\alpha)-\kappa(1-\alpha)}\norm{x(t_0)}E_{\alpha} \left(\frac{\kappa \alpha}{B(\alpha)-\kappa(1-\alpha)}      t^\alpha\right).
$$
}
\end{Proof}
\end{thm}

One crucial aspect of dynamical systems is stability. In control theory, the most popular method to analyze stability is the Lyapunov theory. The work of Lyapunov has had many theoretical implications, generalizations, and practical aspects. In recent years, articles with stability theory of Lyapunov and fractional differential equations have been published. In some papers, for example, in \cite{li2009mittag,li2010stability,sadati2010mittag}, one of the key concepts is the Mittag-Leffler stability that generalizes the usual exponential stability.

\begin{definition}\label{ML_STAB}{\rm
The zero solution of \eqref{system_non} is said to be Mittag-Leffler stable if
\begin{equation}\label{ML-stability}
\norm{x(t)}\leq \left\{ m\left[x(t_0) \right]E_{\alpha}\left(-\lambda t^\alpha\right)   \right\}^{b}
\end{equation}
where $\alpha\in (0,1)$, $\lambda>0, b>0$, $m(0)=0, m(x)\geq 0$ and $m(x)$ is locally Lipschitz function on $x\in\Omega \subset \mathbb{R}^n$ with Lipschitz constant $m_0$.
\hfill{$\Box$}
}\end{definition}

\begin{rem}
Definition \ref{ML_STAB} will be employed in the rest of the chapter to extend some results of stability.  These results are obtained by using the ABC derivative in fractional-order systems, and possible applications in the estimation of state variables will be discussed.
\end{rem}

\section{Lyapunov Direct Method for ABC Fractional-Order Systems}

There are various types of stability of equilibrium points that can be discussed, for example, in \cite{bellman1953stability,rouche1977stability}. One the other hand, if we consider the neighborhood of an equilibrium point in a physical system where the energy is vanishing, then the equilibrium is stable. The Lyapunov's theorems consist of the generalization of this idea because the Lyapunov functions represent a generalization of the energy concept.

Also, these Lyapunov's theorems consist of determining the stability of an equilibrium point by using the properties of the Lyapunov functions indirectly from the knowledge of the solution of the systems \cite{la2012stability}. In this section, we establish the Lyapunov method for fractional-order systems that employ the ABC derivative in their dynamics.

\begin{thm}\label{th_stab1}
Suppose that $a_1, \alpha_2, \alpha_3,a$ and $b$ are arbitrary positive constants, and $x=0$ be an equilibrium point for the system \eqref{system_non} and $\mathcal{D}\subset \mathbb{R}^n$ be a domain containing the origin. Let $V(t,x_t):[0,\infty)\times \mathcal{D}\to \mathbb{R} $ be a continuously differentiable function and locally Lipschitz with respect to $x$ such that
\begin{eqnarray}
\alpha_1 \norm{x}^a \leq V(t,x(t))\leq \alpha_2 \norm{x}^{ab} \label{cond_i} \\
\prescript{ABC}{t_{0}}{\mathcal{D}}^{\alpha}_{t} V(t,x(t))\leq -\alpha_3 \norm{x}^{ab} \label{cond_ii}
\end{eqnarray}
where $t\geq 0$. Then $x=0$ is Mittag-Leffler stable.

\begin{Proof}
{\rm
It follows from inequalities \eqref{cond_i} and \eqref{cond_ii} that
\begin{equation}\label{cond_3}
\prescript{ABC}{t_{0}}{\mathcal{D}}^{\alpha}_{t} V(t,x(t)) +\frac{\alpha_3}{\alpha_2} V(t,x(t))\leq 0
\end{equation}
or
\begin{equation}\label{cond_4}
\prescript{ABC}{t_{0}}{\mathcal{D}}^{\alpha}_{t} V(t,x(t)) +\frac{\alpha_3}{\alpha_2} V(t,x(t))+M(t)=0
\end{equation}
where $M(t)\geq 0$ for $t\geq 0$. Let $\lambda=\dfrac{\alpha_3}{\alpha_2}$, $\gamma=\dfrac{1-\alpha}{B(\alpha)}\lambda$, then inequality \eqref{cond_4} becomes
$$
\prescript{ABC}{t_{0}}{\mathcal{D}}^{\alpha}_{t} V(t,x(t)) + \dfrac{B(\alpha)}{1-\alpha}\gamma V(t,x(t))=-M(t).
$$
Then, taking the Laplace transform to the above equation, we have
\begin{equation*}
 \begin{aligned}
&&\dfrac{B(\alpha)}{1-\alpha}\left[\frac{s^\alpha V(s)-s^{\alpha-1}V(0)+\gamma\left(s^\alpha + \frac{\alpha}{1-\alpha}\right)V(s) }{s^\alpha + \frac{\alpha}{1-\alpha}}   \right]=-M(s) \\
&& \left[s^\alpha\left(\gamma+1\right)+ \frac{\gamma \alpha}{1-\alpha}\right]V(s) =s^{\alpha-1}V(0)-\frac{1-\alpha}{B(\alpha)}s^\alpha M(s)-\frac{\alpha}{B(\alpha)}M(s)
\end{aligned}
\end{equation*}
\begin{equation}\label{final}
V(s)=\kappa_1\dfrac{s^{\alpha-1}}{s^\alpha+\kappa_0}-\left(\kappa_3-\kappa_0\kappa_2\right)\frac{1}{s^\alpha+\kappa_0}M(s)-\kappa_2 M(s)
\end{equation}
where
$$
\kappa_0=\dfrac{\gamma \alpha}{(\gamma+1)(1-\alpha)},\kappa_1=\dfrac{V(0)}{\gamma+1}, \kappa_2=\dfrac{1-\alpha}{B(\alpha)(\gamma+1)},\kappa_3=\dfrac{\alpha}{B(\alpha)(\gamma+1)}.
$$
Applying the inverse Laplace transform to \eqref{final} gives
\begin{eqnarray*}
V(t,x(t))&=&\kappa_1 E_{\alpha}\left(-\kappa_0 t^\alpha\right)-\left(\kappa_3-\kappa_0\kappa_2\right)t^{\alpha-1}E_{\alpha,\alpha}\left(-\kappa_0 t^\alpha\right)\ast M(t)-\kappa_2M(t) \\
&=&\kappa_1 E_{\alpha}\left(-\kappa_0 t^\alpha\right)-\kappa_2M(t)\\
&&\qquad \qquad -\left(\kappa_3-\kappa_0\kappa_2\right)\int_{0}^{t}   (t-\tau)^{\alpha-1}E_{\alpha,\alpha}\left(-\kappa_0 (t-\tau)^\alpha\right)  M(\tau)\:d\tau
\end{eqnarray*}
Since both $t^{\alpha-1}$ and $E_{\alpha,\alpha}\left(-\kappa_0 t^\alpha\right)$ are nonnegative functions, and
$$
\kappa_3-\kappa_0\kappa_2=\frac{\alpha}{B(\alpha)(\gamma+1)^2}>0,
$$
we conclude that
\begin{equation}\label{V_theor}
V(t,x(t))\leq  \kappa_1 E_{\alpha}\left(-\kappa_0 t^\alpha\right)
\end{equation}
From the assumption \eqref{cond_i}, we have the following relationship
\begin{equation}\label{bound_th1}
\norm{x(t)}\leq\left[ \frac{\alpha_2}{\alpha_1 \left(\gamma+1\right)}\norm{x(t_0)}^{ab}E_{\alpha}\left(-\dfrac{\gamma \alpha}{(\gamma+1)(1-\alpha)} t^\alpha\right)\right]^{1/a}
\end{equation}
From which it follows that $x=0$ is Mittag-Leffler stable.
}
\end{Proof}
\end{thm}

\begin{thm}\label{th_stab2}
If the assumptions in Theorem \ref{th_stab1} are satisfied except replacing $\prescript{ABC}{t_{0}}{\mathcal{D}}^{\alpha}_{t}$ by $\prescript{ABR}{t_{0}}{\mathcal{D}}^{\alpha}_{t}$, then one has
\begin{equation}\label{th_stab2_ineq}
\norm{x(t)}\leq\left[ \frac{\alpha_2}{\alpha_1 \left(\gamma+1\right)}\norm{x(t_0)}^{ab}E_{\alpha}\left(-\dfrac{\gamma \alpha}{(\gamma+1)(1-\alpha)} t^\alpha\right)\right]^{1/a}
\end{equation}
where $t\geq t_0$.

\begin{Proof}{\rm
By using Lemma \ref{lema_ineq} and $V(t,x_t)\geq 0$ we obtain $\prescript{ABC}{t_{0}}{\mathcal{D}}^{\alpha}_{t} V(t,x_t) \leq \prescript{ABR}{t_{0}}{\mathcal{D}}^{\alpha}_{t} V(t,x_t)
$ which implies
$$
\prescript{ABC}{t_{0}}{\mathcal{D}}^{\alpha}_{t} V(t,x_t) \leq \prescript{ABR}{t_{0}}{\mathcal{D}}^{\alpha}_{t} V(t,x_t)\leq  -\alpha_3 \norm{x}^{ab}.
$$
Following the proof of Theorem \ref{th_stab1} yields
$$
\norm{x(t)}\leq\left[ \frac{\alpha_2}{\alpha_1 \left(\gamma+1\right)}\norm{x(t_0)}^{ab}E_{\alpha}\left(-\dfrac{\gamma \alpha}{(\gamma+1)(1-\alpha)} t^\alpha\right)\right]^{1/a}
$$
}
\end{Proof}
\end{thm}

This section concludes with an interesting and famous result of differential equations employed in the analysis of differential inequalities and some extensions of the stability results for class-${\mathcal{K}}$ functions.

\begin{definition}\label{class_k}\cite{khalil2002nonlinear} {\rm
A continuous function $\gamma: [0,t)\to [0,\infty)$ is said to belong to class-$\mathcal{K}$ if it is strictly increasing and $\gamma(0)=0$.
\hfill{$\Box$}
}\end{definition}

\begin{lemma}\label{comparison} (The comparison lemma).
Let $0<\alpha<1$, $x(0)=y(0)$ and $$\displaystyle\prescript{ABC}{t_{0}}{\mathcal{D}}^{\alpha}_{t} x(t)\geq \prescript{ABC}{t_{0}}{\mathcal{D}}^{\alpha}_{t} y(t).$$ Then $x(t)\geq y(t)$.

\begin{Proof}
{\rm
It follows from the assumption $\prescript{ABC}{t_{0}}{\mathcal{D}}^{\alpha}_{t} x(t)\geq \prescript{ABC}{t_{0}}{\mathcal{D}}^{\alpha}_{t} y(t)$ that
\begin{equation}\label{eq_ooss}
\prescript{ABC}{t_{0}}{\mathcal{D}}^{\alpha}_{t} x(t)= \prescript{ABC}{t_{0}}{\mathcal{D}}^{\alpha}_{t} y(t)+m(t),
\end{equation}
where $m(t)\geq 0$ for $t\geq 0$. Taking the Laplace transform to \eqref{eq_ooss} yields
\begin{eqnarray*}
\frac{B(\alpha)}{1-\alpha}\dfrac{s^{\alpha-1}}{s^{\alpha}+\dfrac{\alpha}{1-\alpha}}\left[sX(s)-x(0)\right]&=&\frac{B(\alpha)}{1-\alpha}\dfrac{s^{\alpha-1}}{s^{\alpha}+\dfrac{\alpha}{1-\alpha}}\left[sY(s)-y(0)\right]+M(s).
\end{eqnarray*}

Since $x(0)=y(0)$, it follows from the previous equality
$$
X(s)=Y(s)+\frac{1-\alpha}{B(\alpha)}M(s)+\frac{\alpha}{B(\alpha)}  \frac{M(s)}{s^\alpha}
$$
Applying the inverse Laplace transform, one has
$$
x(t)=y(t)+\frac{1-\alpha}{B(\alpha)}m(t)+\frac{\alpha}{B(\alpha)}  \prescript{RL}{t_{0}}{I}^{\alpha}_{t} m(t).
$$
Finally considering that $m(t)\geq 0$ we have that $x(t)\geq y(t)$.
}
\end{Proof}
\end{lemma}

\begin{thm}
Let $x=0$ be an equilibrium point for the nonautonomous fractional-order system \eqref{system_non}. Assume that there exists a Lyapunov function $V(t,x_t)$ and  class-$\mathcal{K}$ functions $\gamma_i$ ($i=1, 2, 3$) such that
\begin{eqnarray}
\gamma_1\left(\norm{x}\right) \leq V(t,x(t))\leq \gamma_2\left(\norm{x}\right) \label{cond_ib} \\
\prescript{ABC}{t_{0}}{\mathcal{D}}^{\alpha}_{t} V(t,x(t))\leq - \gamma_3\left(\norm{x}\right)\label{cond_iib}
\end{eqnarray}
where $\alpha\in(0,1)$. Then, the origin $x=0$ of system \eqref{system_non} is asymptotically stable.

\begin{Proof}{\rm
The proof of this theorem follows the same outline than in the proof of Theorem 11 in \cite{li2009mittag} by using our Lemma \ref{comparison}.
}
\end{Proof}
\end{thm}

\section{Some Useful Inequalities}

Theorem \ref{th_stab1} determines the conditions to conclude Mittag-Leffler stability; however, without a correct Lyapunov function, then we can not determine more about the behavior of the systems.

This section presents some significant inequalities useful in the stability analysis. This analysis employs Lyapunov functions in fractional-order systems whose dynamics are modeled by derivatives with Mittag-Leffler kernel.

\begin{lemma}
Let $x(t)\in \mathbb{R}$ be a differentiable function on $[t_0,T]$ such that $x^\prime(t)\in {L}^{1}(t_0,T)$. Then for any time instant $t\geq t_0$
\begin{equation}\label{ineq_1}
\dfrac{1}{2}   \prescript{ABC}{t_{0}}{\mathcal{D}}^{\alpha}_{t}  x^2(t) \leq x(t)  \prescript{ABC}{t_{0}}{\mathcal{D}}^{\alpha}_{t}  x(t)
\end{equation}
for $0<\alpha<1$ and $t_0<t<T$.

\begin{Proof}
{\rm
To prove relationship \eqref{ineq_1} is equivalent to prove that
\begin{equation}\label{ineq_1b}
 x(t)  \prescript{ABC}{t_{0}}{\mathcal{D}}^{\alpha}_{t}  x(t)-\dfrac{1}{2}   \prescript{ABC}{t_{0}}{\mathcal{D}}^{\alpha}_{t}  x^2(t) \geq 0,\:\forall \alpha \in(0,1).
\end{equation}
Using Definition \ref{1def:ABC}, one has that
\begin{equation}
 x(t)\prescript{ABC}{t_{0}}{\mathcal{D}}^{\alpha}_{t}  x(t)=\dfrac{B(\alpha)}{1-\alpha}\int_{t_0}^{t} x(t)\dot{x}(\tau) E_{\alpha}\left( \dfrac{-\alpha}{1-\alpha}\left(t-\tau\right)^{\alpha} \right)\:d\tau
\end{equation}
and
\begin{equation}
\dfrac{1}{2} \prescript{ABC}{t_{0}}{\mathcal{D}}^{\alpha}_{t}  x^2(t)=\dfrac{B(\alpha)}{1-\alpha}\int_{t_0}^{t}x(\tau)\dot{x}(\tau) E_{\alpha}\left( \dfrac{-\alpha}{1-\alpha}\left(t-\tau\right)^{\alpha} \right)\:d\tau
\end{equation}
Thus, relationship \eqref{ineq_1b} can be written as
\begin{equation}\label{ineq_1c}
\dfrac{B(\alpha)}{1-\alpha}\int_{t_0}^{t}\left[x(t)-x(\tau) \right]\dot{x}(\tau)E_{\alpha}\left( \dfrac{-\alpha}{1-\alpha}\left(t-\tau\right)^{\alpha} \right)\:d\tau\geq 0
\end{equation}
Let $z(\tau)=x(t)-x(\tau)$, then
$$
z^\prime(\tau)=\frac{d z(\tau)}{d\tau}=-\frac{d x(\tau)}{d\tau}.
$$
By using this new variable, the inequality  \eqref{ineq_1c} is rewritten as
\begin{equation}\label{ineq_1d}
\dfrac{B(\alpha)}{1-\alpha}\int_{t_0}^{t} z(\tau)z^\prime(\tau) E_{\alpha}\left( \dfrac{-\alpha}{1-\alpha}\left(t-\tau\right)^{\alpha} \right)\:d\tau \leq 0
\end{equation}
The last inequality can be integrated by parts choosing $dv=z(\tau)z^\prime(\tau)d\tau$, $v=\dfrac{1}{2}z^2(\tau)$ and $u=\dfrac{B(\alpha)}{1-\alpha}E_{\alpha}\left( \dfrac{-\alpha}{1-\alpha}\left(t-\tau\right)^{\alpha} \right)$. According to \cite{martinez2019high}, the derivative of $u$ is given by
$$
\frac{du(\tau)}{d\tau}=\dfrac{\alpha B(\alpha)}{(1-\alpha)^2}\left(t-\tau\right)^{\alpha-1} E_{\alpha,\alpha}\left(\dfrac{-\alpha}{1-\alpha}\left(t-\tau\right)^{\alpha}\right)
$$
Then the integral in the inequality \eqref{ineq_1d} is calculated as follows
\begin{equation*}
\begin{aligned}
&\left. \frac{1}{2} z^2(\tau)\dfrac{B(\alpha)}{1-\alpha}E_{\alpha}\left( \dfrac{-\alpha}{1-\alpha}\left(t-\tau\right)^{\alpha} \right)\right\vert_{\tau=t}-\left. \frac{1}{2} z^2(\tau)\dfrac{B(\alpha)}{1-\alpha}E_{\alpha}\left( \dfrac{-\alpha}{1-\alpha}\left(t-\tau\right)^{\alpha} \right)\right\vert_{\tau=t_0}\\
&-\dfrac{\alpha B(\alpha)}{2(1-\alpha)^2}\int_{t_0}^{t}\left(t-\tau\right)^{\alpha-1}z^2(\tau) E_{\alpha,\alpha}\left(\dfrac{-\alpha}{1-\alpha}\left(t-\tau\right)^{\alpha}\right)\:d\tau
\end{aligned}
\end{equation*}
Since $z(\tau)=x(t)-x(\tau)$, then $z^2(\tau)=0$ if $\tau=t$, and $z^2(t_0)=[x(t)-x(t_0)]^2$. Therefore, we have that
\begin{equation}\label{ineq_1e}
\begin{aligned}
&-\frac{1}{2}[x(t)-x(t_0)]^2\dfrac{B(\alpha)}{1-\alpha}E_{\alpha}\left( \dfrac{-\alpha}{1-\alpha}\left(t-t_0\right)^{\alpha} \right) \\
&-\dfrac{\alpha B(\alpha)}{2(1-\alpha)^2}\int_{t_0}^{t}\left(t-\tau\right)^{\alpha-1} [x(t)-x(t_0)]^2E_{\alpha,\alpha}\left(\dfrac{-\alpha}{1-\alpha}\left(t-\tau\right)^{\alpha}\right)\:d\tau\leq 0
\end{aligned}
\end{equation}
Therefore, in virtue of the positivity of the Mittag-Leffler function for $0<\alpha<1$, then the expression \eqref{ineq_1e} is true and this concludes the proof.
}
\end{Proof}
\end{lemma}

Now, considering the proof in the previous Lemma and properties of the ABC derivative \eqref{ABC_deriv}, we can establish some useful Lemmas. The proofs of these inequalities follow the same outline in \cite{lenka2019time}.

\begin{lemma}
Let $\phi:[t_0,\infty)\to\mathbb{R}$ be a monotonically decreasing and continuously differentiable function. If $x :[t_0,\infty)\to\mathbb{R}$ is a non-negative and continuously differentiable function, then for any time instant $t\geq t_0$ and $0<\alpha<1$:
\begin{equation}\label{lemma_ineq_2}
\prescript{ABC}{t_{0}}{\mathcal{D}}^{\alpha}_{t} \left\{ \phi(t)x(t)\right\}\leq \phi(t) \prescript{ABC}{t_{0}}{\mathcal{D}}^{\alpha}_{t} x(t)
\end{equation}
\end{lemma}

\begin{lemma}
Let $\phi:[t_0,\infty)\to\mathbb{R}$ be a monotonically increasing and continuously differentiable function. If $x :[t_0,\infty)\to\mathbb{R}$ is a non-negative and continuously differentiable function, then for any time instant $t\geq t_0$ and $0<\alpha<1$:
\begin{equation}\label{lemma_ineq_3}
\prescript{ABC}{t_{0}}{\mathcal{D}}^{\alpha}_{t} \left\{ \phi(t)x(t)\right\}\geq \phi(t) \prescript{ABC}{t_{0}}{\mathcal{D}}^{\alpha}_{t} x(t)
\end{equation}
\end{lemma}

\begin{lemma}
Let $\phi:[t_0,\infty)\to\mathbb{R}^n$ be a monotonically decreasing and continuously differentiable vector function. If $x :[t_0,\infty)\to\mathbb{R}^n$ is a non-negative and continuously differentiable vector function, then for any time instant $t\geq t_0$ and $0<\alpha<1$:
\begin{equation}\label{lemma_ineq_4}
\prescript{ABC}{t_{0}}{\mathcal{D}}^{\alpha}_{t} \left\{ \phi^{\intercal}(t)x(t)\right\}\leq \phi^{\intercal}(t) \prescript{ABC}{t_{0}}{\mathcal{D}}^{\alpha}_{t} x(t)
\end{equation}
\end{lemma}

\begin{lemma}
Let $\phi:[t_0,\infty)\to\mathbb{R}$ be a non-negative, monotonically decreasing and continuously differentiable function. If $x :[t_0,\infty)\to\mathbb{R}$ is a continuously differentiable function, then for any time instant $t\geq t_0$ and $0<\alpha<1$:
\begin{equation}\label{lemma_ineq_5}
\prescript{ABC}{t_{0}}{\mathcal{D}}^{\alpha}_{t} \left\{ \phi(t)x^2(t)\right\}\leq 2\phi(t)x(t) \prescript{ABC}{t_{0}}{\mathcal{D}}^{\alpha}_{t} x(t)
\end{equation}
\end{lemma}

\begin{lemma}
Let $\phi:[t_0,\infty)\to\mathbb{R}$ be a non-negative, monotonically decreasing and continuously differentiable function. If $x :[t_0,\infty)\to\mathbb{R}^n$ is a continuously differentiable vector function. Then for any time instant $t\geq t_0$ and $0<\alpha<1$:
\begin{equation}\label{lemma_ineq_6}
\prescript{ABC}{t_{0}}{\mathcal{D}}^{\alpha}_{t} \left\{ \phi(t)x^{\intercal}(t)P x(t)\right\}\leq 2\phi(t)x^{\intercal}(t)P \prescript{ABC}{t_{0}}{\mathcal{D}}^{\alpha}_{t} x(t)
\end{equation}
where $P\in\mathbb{R}^{n\times n}$ is a constant, symmetric, and positive definite matrix.
\end{lemma}

\begin{lemma}
Let $\phi:[t_0,\infty)\to\mathbb{R}^n$ be a vector of differentiable function. Then for any time instant $t\geq t_0$ and $0<\alpha<1$:
\begin{equation}\label{lemma_ineq_7}
\prescript{ABC}{t_{0}}{\mathcal{D}}^{\alpha}_{t} \left\{x^{\intercal}(t)P x(t)\right\}\leq 2 x^{\intercal}(t)P \prescript{ABC}{t_{0}}{\mathcal{D}}^{\alpha}_{t} x(t)
\end{equation}
where $P\in\mathbb{R}^{n\times n}$ is a constant, symmetric, and positive definite matrix.
\end{lemma}

\section{Mittag-Leffler Observers}

Consider the commensurate fractional-order nonlinear systems with a single output described by:
\begin{equation}\label{s1_MLO}
\Sigma_\textbf{x}= \left\{
\begin{aligned}
\prescript{ABC}{t_{0}}{\mathcal{D}}^{\alpha}_{t} \textbf{x}&=f(\textbf{x},\textbf{u})\:,\:\textbf{x}(0)=\textbf{x}_0\\
y&=h(\textbf{x})
\end{aligned}
\right.
\end{equation}
where $0<\alpha<1$, $\textbf{x}\in \mathbb{R}^n$ the state vector; $y\in \mathbb{R}$ the output; $\textbf{u}\in \mathbb{R}^m$ the control input; $f(\textbf{x},\textbf{u}):\mathbb{R}^n\times\mathbb{R}^m\rightarrow \mathbb{R}^n$ a vector function locally Lipschitz in $\textbf{x}$ and uniformly bounded in $\textbf{u}$.

Now, one interesting problem to solve is the estimation of the $n-1$ unknown state variables in the system \eqref{s1_MLO} through another fractional-order dynamic system $\prescript{ABC}{t_{0}}{\mathcal{D}}^{\alpha}_{t} \hat{\textbf{x}}=G\left(\hat{\textbf{x}},y,\textbf{u} \right)$, commonly called observer (or estimator).

Two important characteristics that this algorithm must satisfy are the following \cite{vander}:
\begin{enumerate}
\item If $\hat{\textbf{x}}(t_0)=\textbf{x}(t_0)$ for $t_0>0$, then $\hat{\textbf{x}}(t)=\textbf{x}(t),\: \forall t\geq t_0$.
\item $\hat{\textbf{x}}$ converges to $\textbf{x}$ for $t\to \infty$ sufficiently fast, for any initial conditions $\textbf{x}(t_0)$ and $\hat{\textbf{x}}(t_0)$.
\end{enumerate}

In \cite{martinez2018novel,martinez2019high}, a new class of observers for fractional-order systems is introduced. These observers, called Mittag-Leffler observers, satisfy the conditions 1 and 2 previously mentioned.

\begin{definition}\label{def_observer_MLO}{\rm
Let a fractional dynamical system be described by
\begin{equation}\label{exp_eq}
\prescript{ABC}{t_{0}}{\mathcal{D}}^{\alpha}_{t} \hat{\textbf{x}}=G\left(\hat{\textbf{x}},y,\textbf{u} \right)
\end{equation}
with $\hat{\textbf{x}}\in\mathbb{R}^n$ and $G:\mathbb{R}^n\times \mathbb{R}^p \times \mathbb{R}^m \rightarrow \mathbb{R}^n $ a function continuously differentiable. If there exists an open neighborhood $\mathcal{U}\subset \mathbb{R}^n$ of the origin with $\textbf{x}_0-\hat{\textbf{x}}_0\in \mathcal{U}$ such that
\begin{equation}\label{MLOBS_DEF}
\norm{\textbf{x}-\hat{\textbf{x}}}\leq M \left[E_{\alpha}\left(-\hat{c}t^\alpha \right)\right]^d
\end{equation}
for some positive constants $M, \hat{c}, d$, then the system \eqref{exp_eq} is called a (local) Mittag-Leffler  observer for system \eqref{s1_MLO}.
\hfill{$\Box$}
}\end{definition}

Some comments can be made regarding Definition \ref{def_observer_MLO}:

\begin{enumerate}
\item From the properties of the Mittag-Leffler function, a Mittag-Leffler observer is an asymptotic estimator, i.e.,
$$
\norm{\textbf{x}-\hat{\textbf{x}}}\to 0\:\:\text{as}\:\:t\to\infty
$$
\item If the constant $M$ depends on $\textbf{x}_{0}-\hat{\textbf{x}}_0$, then the inequality \eqref{MLOBS_DEF} is written as:
\begin{equation}\label{MLOBS_DEF_b}
\norm{\textbf{x}-\hat{\textbf{x}}}\leq M_1 \norm{\textbf{x}_{0}-\hat{\textbf{x}}_0}\left[E_{\alpha}\left(-\hat{c}t^\alpha \right)\right]^d
\end{equation}
where $M_1$ is a positive constant.
\item If the open set $\mathcal{U}$ is chosen as the whole space $\mathbb{R}^n$, then \eqref{exp_eq} is called a global Mittag-Leffler observer.
\end{enumerate}

\subsection{Observer for Linear Systems}

Consider the $n$-dimensional equation
\begin{equation}\label{sys_observer_1}
\Sigma_\textbf{x}= \left\{
\begin{aligned}
\prescript{ABC}{t_{0}}{\mathcal{D}}^{\alpha}_{t} \textbf{x}&= A\textbf{x}+B\textbf{u}\\
y&=C\textbf{x}
\end{aligned}
\right.
\end{equation}
where $0<\alpha<1$, $A, B, C$ are given and the input $u(t)$ and the output $y(t)$ are available. However, if the variable $\textbf{x}$ is not available, then the problem is to estimate $\textbf{x}$ from $\textbf{u}$ and $y$ with the knowledge of $A, B, C$. If we know $A$ and $B$ then we can duplicate the original system adding a correcting term as follows \cite{luenberger1964observing}:
\begin{equation}\label{observer_1}
\prescript{ABC}{t_{0}}{\mathcal{D}}^{\alpha}_{t} \hat{\textbf{x}}= A\textbf{x}+B\textbf{u}+\textbf{K}\left(y-C\hat{\textbf{x}}\right).
\end{equation}

Let us define the error between the actual state (system \eqref{sys_observer_1}) and estimated state obtained from system \eqref{observer_1} as follows:
\begin{equation}\label{error_obs1}
\textbf{e}(t)=\textbf{x}(t)-\hat{\textbf{x}}(t)
\end{equation}

Applying the operator ABC \eqref{ABC_deriv} to \eqref{error_obs1} and then substituting \eqref{sys_observer_1} and \eqref{observer_1} into it, we obtain
\begin{equation}\label{error_obs1_dyn}
\prescript{ABC}{t_{0}}{\mathcal{D}}^{\alpha}_{t} \textbf{e}=\left(A-\textbf{K}C\right)\textbf{e}
\end{equation}

Now, let $V=\textbf{e}^\intercal P\textbf{e}$, with $P>0$. Using the inequality \eqref{lemma_ineq_7} and substituting the system \eqref{error_obs1_dyn}, one has that
\begin{eqnarray*}
\prescript{ABC}{t_{0}}{\mathcal{D}}^{\alpha}_{t} V&\leq & 2\textbf{e}^\intercal P \prescript{ABC}{t_{0}}{\mathcal{D}}^{\alpha}_{t}  \textbf{e}\\
&=&2\textbf{e}^\intercal P \left(A-\textbf{K}C\right)\textbf{e}
\end{eqnarray*}
Considering the Lyapunov equation $P \bar{A}+\bar{A}^\intercal P=-Q$, where $\bar{A}=A-\textbf{K}C$, then:
\begin{equation}
\prescript{ABC}{t_{0}}{\mathcal{D}}^{\alpha}_{t} V\leq -\lambda_{\min}( Q )\norm{\textbf{e}}^2
\end{equation}
Therefore, from Theorem \eqref{th_stab1}, the origin $\textbf{e}=0$ is Mittag-Leffler stable and asymptotically stable. In addition, a bound for $\norm{\textbf{e}}$ similar to \eqref{bound_th1} can be obtained.

\begin{rem}
The stability of $\textbf{e}=0$ in system  \eqref{error_obs1_dyn} can be determined by the eigenvalues of the matrix $\bar{A}$. To determine this, we can use similar reasoning as in \cite{matignon1996stability,li2019stability} to establish a region in the complex plane where the poles of $\bar{A}$ can be located.
\end{rem}

\subsection{A High-Gain Nonlinear Observer}

Now, consider the system
\begin{equation}\label{sys_observer_2}
\Sigma_\textbf{x}= \left\{
\begin{aligned}
\prescript{ABC}{t_{0}}{\mathcal{D}}^{\alpha}_{t} \textbf{x}&= A\textbf{x}+\varphi\left(\textbf{x},\textbf{u}\right)  \\
y&=C\textbf{x}
\end{aligned}
\right.
\end{equation}
where
\begin{equation}
  A=\begin{bmatrix}
    0&1&\ldots&0\\
    \vdots&\vdots&\ddots&\vdots\\
    0&0&\ldots&1\\
    0&\ldots&\ldots&0
    \end{bmatrix},\qquad C=\begin{bmatrix}
                              1&0&\ldots&0
                           \end{bmatrix}
\end{equation}
\begin{equation}\label{vector_phi}
 \varphi\left(\textbf{x},\textbf{u}\right)=
 \begin{bmatrix}
 \varphi_1\left(x_1,\textbf{u}\right)\\
 \vdots\\
  \varphi_k\left(x_1,\ldots,x_j,\textbf{u}\right)\\
  \vdots\\
   \varphi_n\left(\textbf{x},\textbf{u}\right)
\end{bmatrix}
\end{equation}

If we assume that the function $\varphi\left(\textbf{x},\textbf{u}\right)$ is Lipschitz with respect to $\textbf{x}$, uniformly with respect to $\textbf{u}$, then together with  the canonical form \eqref{sys_observer_2} we can design many algorithms to solve the problem of estimation in the system \eqref{sys_observer_2}. For example, consider the system \cite{besanccon2007nonlinear}:
\begin{equation}\label{observer_2}
\prescript{ABC}{t_{0}}{\mathcal{D}}^{\alpha}_{t} \hat{\textbf{x}}= A\hat{\textbf{x}}+
\varphi\left(\hat{\textbf{x}},\textbf{u}\right)+\Pi^{-1} \textbf{K}\left(y-C\hat{\textbf{x}}\right)
\end{equation}
where $\textbf{K}$ is a matrix such that $A-\textbf{K} C$ is Hurwitz and
\begin{equation}\label{gains}
\Pi(\theta)=\begin{bmatrix}
      \theta&0&\ldots&0 \\
      0&\theta^2&\ldots&0\\
      \vdots&&\ddots& \\
      0&&\ldots&\theta^n
    \end{bmatrix}
\end{equation}
Now, we verify that algorithm \eqref{observer_2} forms a Mittag-Leffler observer for system \eqref{sys_observer_2}. It is easily verify that:
\begin{eqnarray}
\bar{A}=\Pi(\theta) A \Pi^{-1}(\theta)=\theta^{-1}A\\
\bar{C}=C \Pi^{-1}(\theta)=\theta^{-1}C
\end{eqnarray}

Let $\bar{x}=\Pi(\theta) \textbf{x}$, $\widetilde{x}=\Pi(\theta)  \hat{\textbf{x}}$, $\bar{\varepsilon}=\widetilde{x}-\bar{x}$, $\psi(\bar{\varepsilon})=\varphi\left(\widetilde{x},u\right)-\varphi\left(\bar{x},u\right)$. Then
\begin{eqnarray*}
\prescript{ABC}{t_{0}}{\mathcal{D}}^{\alpha}_{t}\bar{x}&=&\theta^{-1} A \bar{x}+\Pi(\theta)\varphi\left(\Pi^{-1}(\theta)\bar{x},u\right)\\
y&=&\theta^{-1} C \bar{x}\\
\prescript{ABC}{t_{0}}{\mathcal{D}}^{\alpha}_{t}\widetilde{x}&=&\theta^{-1} A\widetilde{x}+\Pi(\theta)\varphi\left(\Pi^{-1}(\theta)\widetilde{x},u\right)+\theta^{-1}\textbf{K}C\left(\bar{x}-\widetilde{x}\right)\\
\prescript{ABC}{t_{0}}{\mathcal{D}}^{\alpha}_{t}\bar{\varepsilon}&=&\theta^{-1}\left(A-\textbf{K} C\right)\bar{\varepsilon}+\Pi\left(\theta\right)\psi\left(\Pi^{-1}(\theta)\bar{\varepsilon}\right)
\end{eqnarray*}

Considering the fact that $A-\textbf{K} C$ is Hurwitz, then there exists a solution of the Lyapunov equation:
\begin{equation}\label{lyap_os}
\left(A-\textbf{K} C\right)^{\intercal}Q+Q\left(A-\textbf{K} C\right)=-I
\end{equation}
Using the solution of \eqref{lyap_os}, let $V=\bar{\varepsilon}^{\intercal}Q\bar{\varepsilon}$, and from inequality \eqref{lemma_ineq_7} one has that
\begin{eqnarray*}
\prescript{ABC}{t_{0}}{\mathcal{D}}^{\alpha}_{t} V&\leq & 2\bar{\varepsilon}^{\intercal}Q \prescript{ABC}{t_{0}}{\mathcal{D}}^{\alpha}_{t}\bar{\varepsilon}\\
&=&2\bar{\varepsilon}^{\intercal}Q \left[\theta^{-1}\left(A-\textbf{K} C\right)\bar{\varepsilon}+\Pi\left(\theta\right)\psi\left(\Pi^{-1}(\theta)\bar{\varepsilon}\right)\right]\\
&=&-\theta^{-1}\norm{\bar{\varepsilon}}^2    +2\bar{\varepsilon}^{\intercal}Q\Pi\left(\theta\right)\psi\left(\Pi^{-1}(\theta)\bar{\varepsilon}\right)\\
&\leq &\left(-\frac{\theta^{-1}}{\lambda_{\max}(Q)}+2\frac{\lambda_{\min}(Q)}{\lambda_{\max}(Q)}\kappa_{\varphi}n^2\right) V
\end{eqnarray*}
where $\kappa_{\varphi}$ is the Lipschitz constant of \eqref{vector_phi}. From this analysis we can conclude the Mittag-Leffler stability by using the Theorem \eqref{th_stab1}.

\subsection{Mittag-Leffler Observers for a Class of Nonlinear Fractional-Order Systems}

Consider the system
\begin{equation}\label{sys_observer_3}
\Sigma_\textbf{x}= \left\{
\begin{aligned}
\prescript{ABC}{t_{0}}{\mathcal{D}}^{\alpha}_{t} \textbf{x}&= A\textbf{x}+\varphi  \\
y&=C\textbf{x}
\end{aligned}
\right.
\end{equation}
where $y\in\mathbb{R}^{n}$, $A\in\mathbb{R}^{n\times n}, C\in\mathbb{R}^{m\times n}$, and the vector $\varphi\left(\textbf{x},y,\prescript{ABC}{t_{0}}{\mathcal{D}}^{\alpha}_{t} y   \right)$ satisfies Lipschitz conditions in $\textbf{x}$, $y$, $\prescript{ABC}{t_{0}}{\mathcal{D}}^{\alpha}_{t} y$. This functions has the following form:
\begin{equation}\label{form_varfi}
\varphi=\varphi_1(y)+\left[\nabla \varphi_2(y)\right]\prescript{ABC}{t_{0}}{\mathcal{D}}^{\alpha}_{t} y+\varphi_3(\textbf{x})
\end{equation}
where $\varphi_1, \varphi_2$ and $\varphi_3$ such that \eqref{sys_observer_3} and \eqref{form_varfi} have unique solutions for any initial $\textbf{x}(0)$.

Now, we introduce the following new variable, $z\in\mathbb{R}^{n}$ as follows:
\begin{equation}\label{eq_z_aux}
z=\textbf{x}-\varphi_2(y)
\end{equation}
Then, applying the ABC derivative, one has
\begin{eqnarray}
\nonumber \prescript{ABC}{t_{0}}{\mathcal{D}}^{\alpha}_{t}  z&=&A\left(z+\varphi_2(y)\right)+\varphi_1(y)+\varphi_3(\textbf{x})\\
&=&Az+\varphi_1(y)+A\varphi_2(y)+\varphi_3\left(z+\varphi_2 \right) \label{eq_z_aux2}
\end{eqnarray}
where $\prescript{ABC}{t_{0}}{\mathcal{D}}^{\alpha}_{t} \varphi_2(y)=\left[\nabla \varphi_2(y)\right]\prescript{ABC}{t_{0}}{\mathcal{D}}^{\alpha}_{t} y$. On the other hand, let
\begin{equation}\label{out_aux}
  y_1=C z=y-C\varphi_2(y)
\end{equation}
be the new output of the system  \eqref{eq_z_aux2}. Based on these calculations and using the variables $y$ and $y_1$ as inputs, the following observer for system  \eqref{eq_z_aux2} is proposed as follows \cite{kou1975exponential}:
\begin{equation}\label{obser_3}
 \prescript{ABC}{t_{0}}{\mathcal{D}}^{\alpha}_{t} \hat{\textbf{x}}=A  \hat{\textbf{x}}+\varphi_1(y)+A\varphi_2(y)+\varphi_3( \hat{\textbf{x}}+\varphi_2)+\textbf{K}\left(y_1-C \hat{\textbf{x}}\right)
\end{equation}
where $\textbf{K}$ is such that
\begin{equation}\label{lyapunov_2}
Q\left(A-\textbf{K}C\right)+\left(A-\textbf{K}C\right)^{\intercal}Q=-P
\end{equation}
with $P,Q$ positive definite and symmetric matrices. Now, by using the Lyapunov direct method we will show that algorithm \eqref{obser_3} is Mittag-Leffler stable. Let $V\left(z-\hat{\textbf{x}}\right)=\left(z-\hat{\textbf{x}}\right)^\intercal Q \left(z-\hat{\textbf{x}}\right)$, then
\begin{eqnarray*}
 \prescript{ABC}{t_{0}}{\mathcal{D}}^{\alpha}_{t} V &\leq &2 \left(z-\hat{\textbf{x}}\right)^\intercal Q \prescript{ABC}{t_{0}}{\mathcal{D}}^{\alpha}_{t}\left(z-\hat{\textbf{x}}\right)\\
 &=&2 \left(z-\hat{\textbf{x}}\right)^\intercal Q\left[\left(A-\textbf{K}C\right)\left(z-\hat{\textbf{x}}\right)+\varphi_3(z+\varphi_2)-\varphi_3(\hat{\textbf{x}}+\varphi_2)\right]\\
 &=&2 \left(z-\hat{\textbf{x}}\right)^\intercal Q\left(A-\textbf{K}C\right)\left(z-\hat{\textbf{x}}\right)\\
 &&\qquad\qquad\qquad\qquad +2 \left(z-\hat{\textbf{x}}\right)^\intercal Q\left[\varphi_3(z+\varphi_2)-\varphi_3(\hat{\textbf{x}}+\varphi_2)\right]\\
\end{eqnarray*}
On the other hand, from the mean-value theorem based on the fundamental theorem of integral calculus, note that
$$
\varphi_3(z+\varphi_2)-\varphi_3(\hat{\textbf{x}}+\varphi_2)=\int_{0}^{1}\nabla\varphi_3(\zeta)\left(z-\hat{\textbf{x}}\right)\:ds,\:\zeta=\hat{\textbf{x}}+\varphi_2+s\left(z-\hat{\textbf{x}}\right)
$$
Then
\begin{eqnarray*}
\prescript{ABC}{t_{0}}{\mathcal{D}}^{\alpha}_{t} V &\leq & -\left(z-\hat{\textbf{x}}\right)^\intercal P \left(z-\hat{\textbf{x}}\right)+2 \left(z-\hat{\textbf{x}}\right)^\intercal Q\int_{0}^{1}\nabla\varphi_3(\zeta)\left(z-\hat{\textbf{x}}\right)\:ds\\
&\leq &-\lambda_{\min} (P)\norm{z-\hat{\textbf{x}}}^2+2 \int_{0}^{1}\abs{\left(z-\hat{\textbf{x}}\right)^\intercal Q \nabla\varphi_3(\zeta)\left(z-\hat{\textbf{x}}\right)}\:ds\\
&\leq &-\lambda_{\min} (P)\norm{z-\hat{\textbf{x}}}^2+2 \int_{0}^{1}\norm{z-\hat{\textbf{x}}}^2 \norm{Q} \sup_{\textbf{x}\in\mathbb{R}^n}\nabla\varphi_3(\textbf{x}) \:ds\\
&=&-\left(\lambda_{\min} (P)-2\norm{Q} \sup_{\textbf{x}\in\mathbb{R}^n}\nabla\varphi_3\right)\norm{z-\hat{\textbf{x}}}^2
\end{eqnarray*}

Then, if
\begin{equation}
\sup_{\textbf{x}\in\mathbb{R}^n}\nabla\varphi_3<\frac{1}{2}\frac{\lambda_{\min} (P)}{\norm{Q}}=\frac{1}{2}\frac{\lambda_{\min} (P)}{\lambda_{\max} (Q)},
\end{equation}
according to Theorem \eqref{th_stab1}, there exists a Mittag-Leffler observer for system \eqref{sys_observer_3}-\eqref{form_varfi}, in particular system \eqref{obser_3}.

\subsection{Review of Numerical Aspects}

For simulation purposes, systems \eqref{s1_MLO}, \eqref{sys_observer_1}, \eqref{sys_observer_2}, \eqref{sys_observer_3}, and Mittag-Leffler observers \eqref{observer_1}, \eqref{observer_2}, \eqref{obser_3} can be implemented through numerical methods, for example by using spectral methods to solve systems of fractional differential equations, in particular the Chebyshev polynomials to obtain an operational matrix of fractional integral using the Clenshaw-Curtis formula, and subsequently solving a system of linear algebraic equations  \cite{baleanu2018chebyshev}.

Other approaches to solving fractional differential equations are the numerical method based on the Adams-Bashforth-Moulton scheme, the Lagrange polynomial interpolation or the homotopy perturbation transform \cite{solis2018novel,coronel2020novel,avalos2019dynamics,yepez2018numerical,kumar2019analysis}. In these papers, the authors show how to solve variable-order fractional systems, so that variable-order Mittag-Leffler observers can be proposed intuitively.

Recently, physical applications modeled by Atangana-Baleanu derivatives have been developed. In these papers, numerical aspects are considered by using the techniques previously mentioned. Interesting comparisons with the discretized convolution integral,  Sudumu-Picard rules, a new version of Adams–Bashforth scheme, and the fixed point method are included \cite{panda2020complex,owolabi2019modelling,baleanu2018nonlinear,gomez2018fractional}.

\section{Conclusion}

Within the framework of the fractional differential equations using the Atangana-Baleanu operators, new results regarding the Mittag-Leffler stability and the Lyapunov Direct Method have been proposed. In this chapter, we have considered the generalization of some inequalities for the Atangana-Baleanu Caputo derivative and their applications in the stability analysis of nonlinear fractional-order systems modeled by ABC and ABR fractional derivatives. Potential applications are addressed by the design of Mittag-Leffler observers for linear and nonlinear fractional-order systems. The numerical implementation is proposed by using various methods published in the literature. To the best of knowledge, there are no results reported of applications and design of the estimation problem of the control theory by using systems with Mittag-Leffler kernel. It is worth to mention that chaotic systems are well defined and can be employed to design the first physical applications of estimation theory and the design of observers.

\bibliographystyle{elsarticle-num} 
\bibliography{biblio_Oscar}

\end{document}